\documentstyle[amscd,amssymb,verbatim]{amsart}
\newcommand{\lrar}[1]{\begin{picture}(50,10)(-25,-5)                          
\put(-25,0){\vector(1,0){50}}
\put(0,5){\makebox(0,0)[b]{\mbox{$#1$}}}
\end{picture}}

\newcommand{\ldar}[1]{\begin{picture}(10,50)(-5,-25)
\put(0,25){\vector(0,-1){50}}
\put(5,0){\mbox{$#1$}}
\end{picture}}

\newcommand{\pa}{\partial}

\newcommand{\bH}{{\bf H}}

\newcommand{\coker}{\operatorname{coker}}

\newcommand{\BB}{{\cal B}}

\renewcommand{\Re}{\operatorname{Re}}
\newcommand{\SL}{\operatorname{SL}}

\newcommand{\G}{{\Bbb G}}

\newcommand{\Coh}{\operatorname{Coh}}

\newcommand{\CC}{{\cal C}}

\newcommand{\ga}{\gamma}
\newcommand{\de}{\delta}
\newcommand{\eps}{\epsilon}

\renewcommand{\ker}{\operatorname{ker}}

\numberwithin{equation}{section}

\newtheorem{thm}{Theorem}[section]
\newtheorem{prop}[thm]{Proposition}
\newtheorem{lem}[thm]{Lemma}
\newtheorem{cor}[thm]{Corollary}
\newenvironment{defi}{\vspace{3mm}\noindent
{\bf Definition.}}{\vspace{3mm}}
\newenvironment{rem}{\vspace{3mm}\noindent
{\bf Remark.}}{\vspace{3mm}}

\newcommand{\Pf}{\noindent {\it Proof}}
\newcommand{\id}{\operatorname{id}}

\newcommand{\ov}{\overline}

\renewcommand{\Im}{\operatorname{Im}}

\newcommand{\rk}{\operatorname{rk}}
\newcommand{\ra}{\rightarrow}

\newcommand{\PP}{{\cal P}}

\newcommand{\SS}{{\cal S}}

\newcommand{\OO}{{\cal O}}

\newcommand{\Hom}{\operatorname{Hom}}

\renewcommand{\a}{\alpha}
\renewcommand{\b}{\beta}

\newcommand{\la}{\lambda}
\newcommand{\th}{\theta}
\newcommand{\C}{{\Bbb C}}
\newcommand{\R}{{\Bbb R}}
\newcommand{\Z}{{\Bbb Z}}

\newcommand{\e}{\operatorname{e}}
\newcommand{\wt}{\widetilde}

\newcommand{\sign}{\operatorname{sign}}
\newcommand{\sub}{\subset}
\newcommand{\ed}{\qed\vspace{3mm}}

\title{Categories of holomorphic vector bundles on noncommutative 
two-tori}
\author{A.~Polishchuk and A.~Schwarz}
\thanks{The work of both authors was partially supported by NSF grants}

\begin{document}
\begin{abstract}
In this paper we study the category of standard holomorphic vector
bundles on a noncommutative two-torus.
We construct a functor from the derived category of such bundles to
the derived category of coherent sheaves on an elliptic curve and prove
that it induces an equivalence with the subcategory of stable
objects. By the homological mirror symmetry for elliptic curves
this implies an equivalence between the derived category of holomorphic
bundles on a noncommutative two-torus and the Fukaya category of the
corresponding symplectic (commutative) torus. 
\end{abstract}

\maketitle

\centerline{\sc Introduction}

\bigskip

In this paper we study holomorphic vector bundles on noncommutative
$2$-dimensional torus $T_{\th}$, where $\th\in\R$,
following the approach of
 \cite{Sp}, \cite{S}, and \cite{DS}.
More precisely, we equip $T_{\th}$ with a complex structure
given by some $\tau\in\C\setminus\R$ and study analogues of
$\bar{\pa}$-connections on vector bundles over the complex noncommutative
torus $T_{\th,\tau}$. 
Recall that vector bundles correspond to projective modules over
the algebra $A_{\th}$ of smooth functions over $T_{\th}$.
As in \cite{DS} we restrict our attention to {\it standard}
holomorphic structures on basic $A_{\th}$-modules,
for which there exist compatible connections of constant curvature.
We do not consider here the problem how to characterize such holomorphic bundles
although it seems that similar to the commutative case the answer should involve some
notion of stability (see remarks at the end of \ref{holbasicsec}).
Motivated by the calculations of
\cite{K} we prove that the {\it derived category}
of standard holomorphic vector bundles on $T_{\th,\tau}$ is equivalent to the
full subcategory of stable objects in the derived category $D^b(X_{\tau})$ of 
coherent sheaves on the complex elliptic curve
$X_{\tau}=\C/(\Z+\tau\Z)$.
In particular, the former category does not depend on $\th$. 
Under this equivalence standard holomorphic
bundles on $T_{\th,\tau}$ land into an abelian subcategory
of $D^b(X_{\tau})$ obtained by certain tilting of the category
$\Coh(X_{\tau})$ of coherent sheaves on $X_{\tau}$ depending on $\th$. Roughly speaking,
to get this abelian subcategory of $D^b(X_{\tau})$,
one has to cut $\Coh(X_{\tau})$ into two subcategories generated by stable bundles
of slopes $<\th$ and $>\th$ respectively (we assume that $\th$ irrational) and then
reassemble these subcategories in a different way into a new abelian category.
The general framework for such constructions is provided by torsion theory of \cite{HRS}.
These $t$-structures were considered by Bridgeland in connection with
general stability conditions for derived categories.

In physical terms the objects of our study are D-branes of B-type in
SUSY Yang-Mills theory on noncommutative $2$-tori.
Although the corresponding derived categories
are the same for all $\th$, the physical branes (stable objects) do depend
on $\th$. This resembles the picture 
where stability conditions on the derived category of
coherent sheaves over a Calabi-Yau manifold vary with
K\"ahler structure.
(see \cite{Doug} for more general discussion).

An interesting feature of the above equivalence of categories is that it is compatible with
Morita equivalences and the action of
$\SL_2(\Z)$ on $D^b(X_{\tau})$ defined by Mukai (see \cite{M}).
Namely, for every $g\in\SL_2(\Z)$ one can consider the fractional-linear action of $g$ 
on parameter $\th$. The noncommutative tori $T_{g\th}$ and $T_{\th}$ are Morita equivalent,
i.e., there exists an $A_{g\th}-A_{\th}$-bimodule inducing an equivalence of categories
of modules over $A_{g\th}$ and $A_{\th}$. This equivalence extends to an equivalence between
categories of holomorphic bundles on $T_{g\th,\tau}$ and $T_{\th,\tau}$. 
Now one can compare equivalences between the corresponding derived categories with 
$D^b(X_{\tau})$.
It turns out that they differ by the action of the autoequivalence of $D^b(X_{\tau})$
corresponding to $g^t$. Thus, our picture can be considered as a 
``noncommutative explanation'' of the existence of the Fourier-Mukai transform.

We consider two approaches to constructing the above equivalences.
One uses the explicit calculation of structure constants
of products of noncommutative theta functions in \cite{DS}.
Another is based on a construction of the Fourier-Mukai type functor
from the category of holomorphic vector bundles on $T_{\th,\tau}$
to the category of complexes of sheaves of $\OO$-modules
with coherent cohomology on $X_{\tau}$.

We restrict ourselves to considering only standard holomorphic structures on basic
projective modules over $A_{\th}$, which corresponds to considering only stable
vector bundles on the commutative elliptic curve (or more generally, simple coherent
sheaves). However, we believe that the functor we construct 
extends to all holomorphic bundles on $T_{\th}$.
It is also clear that some of the results of this paper
can be generalized to higher-dimensional noncommutative tori. Note that
in this generalization sometimes one has to replace the category of 
coherent sheaves on the corresponding commutative complex torus $X$
by the twisted category associated with
some cohomology class in $H^2(X,\OO^*)$ 
(this twisting is always trivial in $2$-dimensional case).

In the paper \cite{K} that motivated us
the category of standard holomorphic bundles on $T_{\th}$
is compared with the Fukaya category of a symplectic torus. In view of
the equivalence we establish in this paper, such a comparison essentially
reduces to the usual homological mirror symmetry for an elliptic curve
proved in \cite{PZ}. It is remarkable that the functor from the 
Fukaya category to the category of
holomorphic bundles on $T_{\th}$ considered in \cite{K} seems to be simplier
than the similar functor in the commutative case. It would be
interesting to study the relation between the higher-dimensional
generalization of this functor and Fukaya's work on noncommutative
mirror symmetry \cite{F}.

We should mention that a different construction of
categories associated with noncommutative elliptic curves was
proposed by Soibelman and Vologodsky in \cite{SV}. The categories
$\BB_q$
they construct depend only on one parameter $q=\exp(2\pi i\th)$,
while our categories depend on two parameters $\th$ and $\tau$.
It seems plausible that the categories $\BB_q$ can be obtained
as some asymptotic limits of our categories when $\tau$ goes to infinity.
On the other hand, the categories $\BB_q$ are in some sense 
degenerations of the categories
of coherent sheaves on commutative elliptic curves $\C^*/q^{\Z}$ (where
$|q|<1$) as $|q|\to 1$. The general philosophy of this degeneration
procedure was outlined in \cite{So}.

It would be also interesting to understand the role of holomorphic structures
in Manin's program that is supposed to relate the real multiplication for
noncommutative two-tori with arithmetics of real quadratic fields
(see \cite{Manin}).

\noindent
{\it Acknowledgment.} We are grateful to Marc Rieffel, Yan Soibelman and Mauro Spera for
useful discussions. Part of this research was carried out during the
stay of both authors at the Institut des Hautes Etudes Scientifiques. 
We thank this institution for the hospitality and for stimulating working conditions.

\section{Preliminaries}

\subsection{Basic modules}\label{basicsec}

  Let $\th$ be a real number, $T_{\th}$ be the corresponding
$2$-dimensional noncommutative torus.
By the definition, the algebra $A_{\th}$ of smooth functions on $T_{\th}$ 
consists of formal linear combinations
$\sum_{(n_1,n_2)\in\Z^2} a_{n_1,n_2}U_1^{n_1}U_2^{n_2}$ with
the coefficient function $(n_1,n_2)\mapsto a_{n_1,n_2}$ rapidly
decreasing at infinity. The product is defined using the rule
$$U_1U_2=\exp(2\pi i\th)U_2U_1.$$

We are interested in finitely generated projective 
right $A_{\th}$-modules. It is known 
(see \cite{Co}, \cite{R-ct}) that every such module
is isomorphic to one of the following modules $E_{n,m}(\th)$,
where $(m,n)$ is a pair of integers such that $n+m\th\neq 0$.
Assume first that $m\neq 0$. Then $E_{n,m}(\th)$ is defined as
the Schwartz space $\SS(\R\times\Z/m\Z)$ equipped with the following
right action of $A_{\th}$:
$$fU_1(x,\a)=f(x-\frac{n+m\th}{m},\a-1),\
fU_2(x,\a)=\exp(2\pi i (x-\frac{\a n}{m}))f(x,\a),$$
where $x\in\R$, $\a\in\Z/m\Z$.
For $m=0$ and $n\neq 0$ we set $E_{n,0}(\th)=A_{\th}^{|n|}$ with the obvious
right $A_{\th}$-action. Let us denote by
$$\deg(E_{n,m}(\th))=m,$$
$$\rk(E_{n,m}(\th))=n+m\th,$$
$$\mu(E_{n,m}(\th))=\frac{m}{n+m\th}$$
the natural analogues of degree, rank and slope.
%In \cite{} the number $|\rk(E_{n,m}(\th))|$ is
%called the non-commutative dimension of the module $E_{n,m}(\th)$. 
The modules $E_{-n,-m}(\th)$ and $E_{n,m}(\th)$ are equivalent, however, sometimes
we would like to distinguish them by introducing the $\Z/2\Z$-grading. 
Namely, we will say that the module $E$ is even if $\rk E>0$ and that it is odd if
$\rk E<0$. Thus, in the category of $\Z/2\Z$-graded projective $A_{\th}$-modules we have
$$E_{-n,-m}(\th)=\Pi E_{n,m}(\th),$$
where $E\mapsto \Pi E$ denotes the parity change.
In the case when $\th$ is rational, we could also try to use the above
formulas to define modules associated with $(m,n)$ such that $n+m\th=0$; however, these modules
are not projective.

{\it Basic modules} are the modules $E_{n,m}(\th)$ with $m$ and $n$
relatively prime. It is known that for $(m,n)$ relatively prime
one has $E_{nd,md}(\th)\simeq E_{n,m}(\th)^{\oplus d}$.

It is often convenient to extend the pair $(m,n)$ to a matrix in $\SL_2(\Z)$.
When we need to make this choice we will use the following notation:
for a matrix $g=\left(\matrix a & b \\ c & d\endmatrix\right)\in\SL_2(\Z)$ 
we set
$E_g(\th)=E_{d,c}(\th)$. In particular, for $g=1$, we get
$E_1(\th)=E_{1,0}(\th)=A_{\th}$. In accordance with our formulas for
the degree and rank of modules we set for a matrix $g$ as above
$$\deg(g)=c,$$
$$\rk(g,\th)=c\th+d,$$
so that $\deg(E_g(\th))=\deg(g)$ and $\rk(E_g(\th))=\rk(g,\th)$.
Note that the function $g\mapsto \rk(g,\th)$ satisfies the following
cocycle condition
\begin{equation}\label{cocycle}
\rk(g_1g_2,\th)=\rk(g_1,g_2\th)\rk(g_2,\th),
\end{equation}
where we denote by 
$\th\mapsto g\th=(a\th+b)/(c\th+d)$ the natural action of $\SL_2(\Z)$
on the set of $\th$'s (it is always well-defined for irrational $\th$,
and partially defined for rational $\th$). 
We will also need the following identity which can be easily checked:
\begin{equation}\label{degid}
\deg(g_2g_1^{-1})\rk(g_3,\th)-\deg(g_3g_1^{-1})\rk(g_2,\th)+\deg(g_3g_2^{-1})\rk(g_1,\th)=0.
\end{equation}
It is well-known that endomorphisms of the basic right $A_{\th}$-module
$E_g(\th)$ can be identified with the algebra $A_{g\th}$. The
corresponding left action of $A_{g\th}$ on $E_g(\th)$ is given by
the formulas
$$U_1f(x,\a)=f(x-\frac{1}{c},\a-a),\ 
U_2f(x,\a)=\exp(2\pi i(\frac{x}{c\th+d}-\frac{\a}{c}))f(x,\a),$$
where $g=\left(\matrix a & b \\ c & d\endmatrix\right)$. 

\subsection{Pairings}

We have a canonical pairing
\begin{equation}\label{bpair}
b:E_{g^{-1}}(g\th)\otimes_{\C} E_g(\th)\ra\C,
\end{equation}
defined by the formula 
$$b(f_1\otimes f_2)=\sum_{\a\in\Z/\deg(g)\Z}\int_{x\in\R} f_1(\frac{x}{\rk(g,\th)},\a)f_2(x,-a\a)dx,$$
where $f_1\in E_{g^{-1}}(g\th)$, $f_2\in E_g(\th)$, $g$ is given by the matrix above.

\begin{lem} For every $U\in A_{g\th}$, $V\in A_{\th}$ one has
$$b(f_1 U\otimes f_2)=b(f_1\otimes Uf_2),$$
$$b(Vf_1\otimes f_2)=b(f_1\otimes f_2V).$$
\end{lem}

\Pf . It suffices to check these equalities when $U$ (resp. $V$) is one of the two generators of 
the algebra $A_{g\th}$ (resp. $A_{\th}$), in which case it is straighforward.
\ed

%The above notation is especially useful when considering
%tensor products of basic modules as the reader can see from
%the following result.

\begin{prop}\label{pairprop} 
(a) For every $g_1,g_2\in\SL_2(\Z)$ one has a well-defined
pairing of right $A_{\th}$-modules 
$$t_{g_1,g_2}:E_{g_1}(g_2\th)\otimes_{\C} E_{g_2}(\th)\ra E_{g_1g_2}(\th)$$
given by the following formulas. If $g_1\neq 1$, $g_2\neq 1$ and 
$g_1g_2\neq 1$ then
\begin{align*}
&t_{g_1,g_2}(f_1\otimes f_2)(x,\a)=\\
&\sum_{n\in\Z}f_1(\frac{x}{\rk(g_2,\th)}+\frac{\rk(g_1,g_2\th)}{c_1}
\left(\frac{c_2d_{12}\a}{c_{12}}-n\right),a_1d_{12}\a-n)
f_2(x-\frac{d_{12}\a}{c_{12}}+\frac{n}{c_2},a_2n),
\end{align*}
where $x\in\R$, $\a\in \Z/c_{12}\Z$,
$g_i=\left(\matrix a_i & b_i \\ c_i & d_i \endmatrix\right)$,
$g_1g_2=\left(\matrix a_{12} & b_{12} \\ c_{12} & d_{12} \endmatrix\right)$.
The pairing 
$$t_{g,1}:E_g(\th)\otimes_{\C} E_1(\th)\ra E_g(\th)$$
is given by the right action of the algebra $A_{\th}=E_1(\th)$ on 
$E_g(\th)$. The pairing
$$t_{1,g}:E_1(g\th)\otimes_{\C} E_g(\th)\ra E_g(\th)$$
is set to be the left action of the algebra $A_{g\th}=E_1(g\th)$
on $E_g(\th)$ defined before. Finally,
the pairing 
$$t_{g^{-1},g}:E_{g^{-1}}(g\th)\otimes_{\C} E_g(\th)\ra E_1(\th)=A_{\th}$$
is given by the formula
$$t_{g^{-1},g}(f_1\otimes f_2)=\sum_{(n_1,n_2)\in\Z^2}
U_1^{n_1}U_2^{n_2}b(U_2^{-n_2}U_1^{-n_1}f_1\otimes f_2),$$ 
where $b$ is the pairing (\ref{bpair}).

\noindent
(b) For every triple of elements 
$g_1,g_2,g_3\in\SL_2(\Z)$ the following diagram is commutative
\begin{equation}
\begin{array}{ccc}
E_{g_1}(g_2g_3\th)\otimes E_{g_2}(g_3\th)\otimes E_{g_3}(\th) &\lrar{t_{g_1,g_2}\otimes\id}& 
E_{g_1g_2}(g_3\th)\otimes E_{g_3}(\th)\\
\ldar{\id\otimes t_{g_2,g_3}} & &\ldar{t_{g_1g_2,g_3}}\\
E_{g_1}(g_2g_3\th)\otimes E_{g_2g_3}(\th) &\lrar{t_{g_1,g_2g_3}} & E_{g_1g_2g_3}(\th)
\end{array}
\end{equation}

\noindent
(c) The map $E_{g_1}(g_2\th)\ra\Hom_{A_{\th}}(E_{g_2}(\th),E_{g_1g_2}(\th))$ induced by $t_{g_1,g_2}$
is an isomorphism.
\end{prop}

Before starting the proof it is convenient to rewrite the formula for $t_{g_1,g_2}$ a little bit.

\begin{lem}\label{pairlem}
(a) For $g_1,g_2\in\SL_2(\Z)$ such that $g_1\neq 1$, $g_2\neq 1$, $g_1g_2\neq 1$ one has
\begin{align*}
&t_{g_1,g_2}(f_1\otimes f_2)(x,\a)=\sum_{\a_1\in\Z/c_1\Z,\a_2\in\Z/c_2\Z}\\
&\sum_{n\in I_{g_1,g_2}(\a_1,\a_2,\a)}
f_1(\frac{x}{\rk(g_2,\th)}+\frac{\rk(g_1,g_2\th)n}{c_1c_{12}},\a_1)
f_2(x-\frac{n}{c_2c_{12}},\a_2),
\end{align*}
where $g_1$, $g_2$ and $g_1g_2$ are given by the same matrices as in the above proposition,
\begin{align*}
&I_{g_1,g_2}(\a_1,\a_2,\a)=\{n\in\Z|\ n\equiv -c_1\a+c_{12}\a_1 (c_{12}c_1),\\ 
&n\equiv c_2d_{12}\a-c_{12}d_2\a_2 (c_{12}c_2)\}.
\end{align*}

\noindent
(b) For $n\in I_{g_1,g_2}(\a_1,\a_2,\a)$ one has $n\equiv-c_1\a_2+d_1c_2\a_1 (c_1c_2)$.

\noindent
(c) One has $I_{g_2,-g_2^{-1}g_1^{-1}}(\a_1,\a_2,\a_3)=I_{g_1,g_2}(-a_1\a_3,-\a_1,-a_{12}\a_2)$.
\end{lem}

\noindent
{\it Proof of lemma}. (a) In the formula of Proposition \ref{pairprop}(a) one should
change the summation variable by setting $n'=c_2d_{12}\b-c_{12}n$. Then the congruences
$$a_1d_{12}\a-n\equiv\a_1(c_1),\ a_2n\equiv\a_2(c_2)$$
are equivalent to the congruences
$$n'\equiv c_2d_{12}\a-c_{12}d_2\a_2(c_{12}c_2),\ 
n'\equiv (c_2-a_1c_{12})d_{12}\a+c_{12}\a_1 (c_{12}c_1)$$
(we used the fact that $d_2\equiv a_2^{-1}(c_2)$).
It remains to use the relation $(c_2-a_1c_{12})d_{12}=-c_1a_{12}d_{12}\equiv -c_1 (c_{12}c_1)$
to see that these congruences are equivalent to $n'\in I_{g_1,g_2}(\a_1,\a_2,\a)$.

\noindent
(b) Let us pick a representative in $\Z$ for $\a\in\Z/c_{12}\Z$.
Since $c_2d_{12}-c_{12}d_2=-c_1$ we have the congruence
$$n\equiv -c_1\a+c_{12}d_2(\a-\a_2)(c_{12}c_2)$$
for $n\in I_{g_1,g_2}(\a_1,\a_2,a)$.
Hence, $k:=(n+c_1\a)/c_{12}$ is an integer and we have
$$k\equiv d_2(\a-\a_2) (c_2),\  k\equiv\a_1 (c_1).$$
Using these two relations we get
$$n+c_1\a=c_{12}k=c_1a_2k+d_1c_2k\equiv c_1a_2d_2(\a-\a_2)+d_1c_2\a_1 (c_1c_2)\equiv
c_1(\a-\a_2)+d_1c_2\a_1 (c_1c_2)$$
as required.

\noindent
(c) This follows easily from (b).
\ed

\noindent
{\it Proof of Proposition \ref{pairprop}}. 
The pairing $t_{g_1,g_2}$ for $g_1\neq 1$, $g_2\neq 1$, $g_1g_2\neq 1$ is essentially the same
as in section 2 of \cite{DS}. We just replaced left modules by right modules everywhere.
Now the results of \cite{DS} imply that $t_{g_1,g_2}$ for $g_1\neq 1$, $g_2\neq 1$ gives rise
to a map of $A_{g_1g_2\th}$-$A_{\th}$-bimodules (in fact, an isomorphism)
\begin{equation}\label{bimisom}
E_{g_1}(g_2\th)\otimes_{A_{g_2\th}}E_{g_2}(\th)\ra E_{g_1g_2}(\th).
\end{equation}
This proves the associativity in (b) for the case when one of the elements $g_i$ is equal to $1$,
and two other elements $g_j$ and $g_k$ are different from $1$ and are not inverse to each other.
It is also easy to see directly that the pairing $t_{g^{-1},g}$ gives rise to a map of
$A_{\th}$-$A_{\th}$-bimodules (also an isomorphism)
$$E_{g^{-1}}(g\th)\otimes_{A_{g\th}}E_g(\th)\ra A_{\th},$$
which implies the associativity in (b) for the triples $(g^{-1},g,1)$, $(g^{-1},1,g)$ and
$(1,g^{-1},g)$. 

In the case when all the elements $g_1$, $g_2$, $g_3$, $g_1g_2$, $g_2g_3$, $g_1g_2g_3$
are non-trivial the associativity can be checked by a straightforward calculation
using Lemma \ref{pairlem}(a),(b). More precisely, if we denote
$g_i=\left(\matrix a_i & b_i \\ c_i & d_i \endmatrix\right)$,
$g_ig_j=\left(\matrix a_{ij} & b_{ij} \\ c_{ij} & d_{ij} \endmatrix\right)$,
$g_1g_2g_3=\left(\matrix a_{123} & b_{123} \\ c_{123} & d_{123} \endmatrix\right)$, then
the associativity follows from the statement that for 
every collection $(\a_1,\a_2,\a_3,\a)$, where $\a_i\in\Z/c_i\Z$, $\a\in\Z/c_{123}\Z$, the 
map $(m,n)\mapsto ((c_3m+c_{123}n)/c_{23},(c_2m-c_1n)/c_{23})$ gives a bijection
\begin{align*}
&\cup_{\a_{23}\in\Z/c_{23}\Z}I_{g_1,g_2g_3}(\a_1,\a_{23},\a)\times I_{g_2,g_3}(\a_2,\a_3,\a_{23})
\wt{\ra} \\
&\cup_{\a_{12}\in\Z/c_{12}\Z}I_{g_1g_2,g_3}(\a_{12},\a_3,\a)\times I_{g_1,g_2}(\a_1,\a_2,\a_{12}).
\end{align*}

The remaining interesting cases are: (i) $g_1g_2g_3=1$, all elements $g_1$, $g_2$, $g_3$ are non-trivial;
(ii) $g_1g_2=1$, all elements $g_2$, $g_3$, $g_2g_3$ are non-trivial; (ii)' $g_2g_3=1$, all elements
$g_1$, $g_2$, $g_1g_2$ are non-trivial.
The case (i) boils down to the identity
$$b(t_{g_1,g_2}(f_1\otimes f_2)\otimes f_3)=b(f_1,t_{g_2,g_3}(f_2\otimes f_3))$$
which is not difficult to check directly.
Let us prove the identity in the case (ii) (the case (ii)' is similar).
We can assume that $\th$ is irrational (the case of rational $\th$ will follow
by continuity).
Then we can easily show that the associativity holds up to a constant. Indeed,
for fixed $f_1\in E_{g_2^{-1}}(g_2g_3\th)$, $f_2\in E_{g_2}(g_3\th)$, the map
$$f_3\mapsto t_{g_2^{-1},g_2g_3}(f_1\otimes t_{g_2,g_3}(f_2\otimes f_3)$$
is a morphism of right $A_{\th}$-modules $E_{g_3}(\th)\ra E_{g_3}(\th)$,
hence it is given by the left action of some element $t'(f_1,f_2)\in A_{g_3\th}$. 
Moreover, the map $t'$ gives a morphism of $A_{g_3\th}-A_{g_3\th}$-bimodules
$$E_{g_2^{-1}}(g_2g_3\th)\otimes_{A_{g_2g_3\th}}E_{g_2}(g_3\th)\ra A_{g_3\th}.$$
Since $t_{g_2^{-1},g_2}$ induces an isomorphism of the bimodule on the left
with $A_{g_3\th}$ we see that this morphism differs from $t_{g_2^{-1},g_2}$ by
the multiplication with a central element in $A_{g_3\th}$. Since we assumed that $\th$ is
irrational, this central element should be a constant.
Thus, it suffices to prove the following identity:
\begin{equation}\label{assoeq}
t_{g_2^{-1},g_2g_3}(f_1\otimes t_{g_2,g_3}(f_2\otimes f_3)(0,0)=
(t_{g_2^{-1},g_2}(f_1\otimes f_2)f_3)(0,0).
\end{equation}
Also, we can assume that $f_1(x,\a)=0$ for $\a\not\equiv 0(c_1)$ and
denote $f_1(x)=f_1(x,0)$. Then the LHS of (\ref{assoeq}) can be rewritten as
\begin{align*}
&\sum_{m,n\in\Z}f_1(\frac{m}{\rk(g_2,g_3\th)})f_2(m+\frac{n\rk(g_2,g_3\th)}{c_2},n)
f_3(-\frac{n}{c_3},-a_3n)=\\
&\sum_{n\in\Z}(U_1^nf_3)(0,0)\cdot\sum_{m\in\Z}
f_1(\frac{m}{\rk(g_2,g_3\th)})(f_2U_1^{-n})(m,0).
\end{align*}
Now we can rewrite the inner sum using the Poisson summation formula:
\begin{align*}
&\sum_{m\in\Z}f_1(\frac{m}{\rk(g_2,g_3\th)})(f_2U_1^{-n})(m,0)=
\sum_{m\in\Z}\int\exp(2\pi imy)f_1(\frac{y}{\rk(g_2,g_3\th)})(f_2U_1^{-n})(y,0)dy=\\
&\sum_{m\in\Z}\int f_1(\frac{y}{\rk(g_2,g_3\th)})(f_2U_1^{-n}U_2^m)(y,0)dy=
\sum_{m\in\Z} b(f_1,f_2U_1^{-n}U_2^{-m}).
\end{align*}
Therefore, the LHS of (\ref{assoeq}) is equal to
$$\sum_{m,n\in\Z}(U_1^nf_3)(0,0)b(U_1^{-n}U_2^{-m}f_1,f_2)=
\sum_{m,n\in\Z}(U_2^mU_1^nf_3)(0,0)b(U_1^{-n}U_2^{-m}f_1,f_2)$$
which is the definition of the RHS of (\ref{assoeq}).

Part (c) follows from the well-known fact that all basic modules considered as bimodules
are invertible and from the isomorphism (\ref{bimisom}).
\ed

\begin{rem} The cocycle condition \ref{cocycle} implies that the pairing
$t_{g_1,g_2}$ is even (where we define the $\Z/2\Z$-grading
of the tensor product in the usual way).
\end{rem}

\begin{cor}\label{compcor}
(a) For every pair of basic modules $E_g(\th)$, $E_{h}(\th)$ there is a canonical
isomorphism 
\begin{equation}\label{homisom}
\Hom_{A_{\th}}(E_g(\th),E_{h}(\th))\simeq E_{hg^{-1}}(g\th)
\end{equation}
of $A_{h\th}-A_{g\th}$-bimodules. 

\noindent (b)
For a triple of basic modules 
$(E_g(\th), E_{h}(\th), E_{k}(\th))$ the canonical composition map 
$$\Hom_{A_{\th}}(E_{h}(\th),E_{k}(\th))\otimes \Hom_{A_{\th}}(E_g(\th),E_{h}(\th))
\ra \Hom_{A_{\th}}(E_g(\th),E_{k}(\th))$$
is identified via isomorphisms (\ref{homisom}) with the pairing $t_{kh^{-1},hg^{-1}}$.
\end{cor}

\begin{rem}
Recall that the Morita
equivalence associated with a basic module $E_{g_0}(\th)$
is the functor $M\mapsto M\otimes_{A_{g_0\th}}E_{g_0}(\th)$ from the category of 
right $A_{g_0\th}$-modules
to that of right $A_{\th}$-modules. According to (\ref{bimisom}) it sends
$E_g(g_0\th)$ to $E_{gg_0}(\th)$. The corresponding isomorphism
$$\Hom_{A_{g_0\th}}(E_g(g_0\th),E_h(g_0\th))\simeq\Hom_{A_{\th}}(E_{gg_0}(\th),E_{hg_0}(\th))$$
is compatible with the identification of both sides with
$E_{hg^{-1}}(gg_0\th)$ given in the above corollary.
\end{rem}

\section{Holomorphic structures}

\subsection{Standard holomorphic structures on basic modules}\label{holbasicsec}

Let us fix a complex number $\tau$ such that $\Im(\tau)\neq 0$. 
We will think about $\tau$ as a complex structure on
the noncommutative torus $T_{\th}$. 
Namely, $\tau$ defines a one-dimensional subalgebra in the Lie
algebra of derivations of the algebra $A_{\th}$ spanned by the derivation $\de_{\tau}$
given by 
\begin{equation}\label{detau}
\de_{\tau}(\sum_{(n_1,n_2)\in\Z^2}a_{n_1,n_2}U_1^{n_1}U_2^{n_2})=
2\pi i \sum_{(n_1,n_2)\in\Z^2}(n_1\tau+n_2)a_{n_1,n_2}U_1^{n_1}U_2^{n_2}.
\end{equation}
We denote by $T_{\th,\tau}$ the noncommutative torus $T_{\th}$ equipped
with this complex structure.

\begin{defi} (i) A {\it holomorphic structure} on a right $A_{\th}$-module $E$ 
(compatible with the complex structure on $T_{\th,\tau}$) is an operator
$\ov{\nabla}:E\ra E$ satisfying the following Leibnitz identity:
$$\ov{\nabla}(ea)=\ov{\nabla}(a)\cdot a+e\cdot\de_{\tau}(a),$$
where $e\in E$, $a\in A_{\th}$.

\noindent 
(ii) If $E$ and $E'$ are right $A_{\th}$-modules equipped with holomorphic structures, then
we say that a morphism $f:E\ra E'$ of right $A_{\th}$-modules is {\it holomorphic}
if $\ov{\nabla}(fe)=f(\ov{\nabla}e)$.

\noindent
(iii) For a right $A_{\th}$-module $E$ equipped with a holomorphic structure $\ov{\nabla}$ we define
the cohomology as follows: 
$$H^*(E)=H^*(E,\ov{\nabla}):=H^0(E,\ov{\nabla})\oplus H^1(E,\ov{\nabla}),\ \text{where}$$
$$H^0(E,\ov{\nabla}):=\ker(\ov{\nabla}:E\ra E),\ H^1(E,\ov{\nabla}):=\coker(\ov{\nabla}:E\ra E).$$ 
\end{defi}

It is natural to call projective right $A_{\th}$-modules equipped with
holomorphic structures {\it holomorphic bundles} on the
complex noncommutative torus $T_{\th,\tau}$.

We are going to introduce a family of holomorphic structures (depending on one complex
parameter) on every basic module. Namely, for $E=E_{n,m}(\th)$ with $m\neq 0$ we set
$$\ov{\nabla}_z(f)=\frac{\pa f}{\pa x}+2\pi i(\tau\mu(E)x+z)f,$$
where $\mu(E)=m/(n+m\th)$, $f\in E$, $z\in\C$. 
It is easy to see that for every $z$ this is
a holomorphic structure on $E$. For $E=A_{\th}$ we define holomorphic
structures on the trivial module $\ov{\nabla}_z:A_{\th}\ra A_{\th}$ by setting
$$\ov{\nabla}_z(a)=2\pi i z\cdot a+\de_{\tau}(a).$$ 
We will denote by $E_g^z(\th)$ the basic module $E_g$ equipped with the holomorphic structure
$\ov{\nabla}_z$. 

\begin{defi} A {\it standard} holomorphic structure on a basic module $E_g$ is
one of the structures $\ov{\nabla}_z$. A {\it standard} holomorphic bundle on $T_{\th}$ is
one of the bundles $E_g^z(\th)$.
\end{defi}

Note that there is a natural notion of equivalence of holomorphic structures
on the same module $E$: two such structures $\ov{\nabla}$ and $\ov{\nabla}'$ are equivalent
if there exists a holomorphic isomorphism $(E,\ov{\nabla})\ra (E,\ov{\nabla}')$. In other words,
the group of automorphisms of $E$ as a right $A_{\th}$-module acts on the space of holomorphic
structures on $E$, so that equivalence classes are orbits under this action. 

\begin{prop}\label{leftLeiblem} 
(a) For every basic right $A_{\th}$-module $E=E_g(\th)$ and every $z\in\C$ 
the operator $\ov{\nabla}_z:E\ra E$ satisfies the
following Leibnitz rule with respect to the left action of $A_{g\th}$ on $E$:
$$\ov{\nabla}_z(be)=b\cdot\ov{\nabla}_z(e)+\frac{1}{\rk(E)}\de_{\tau}(b)\cdot e,$$
where $b\in A_{g\th}$, $e\in E$.

\noindent
(b) The holomorphic structures
$\ov{\nabla}_{z'}$ and $\ov{\nabla}_z$ are equivalent
if and only if $z'\equiv z\mod\frac{1}{\rk(E)}(\Z+\tau\Z)$.

\noindent
(c) Let $E_0=E_g(\th)$ be a basic right $A_{\th}$-module equipped with a standard
holomorphic structure $\ov{\nabla}_z$, $E$ be a right $A_{g\th}$-module. Then 
for every holomorphic structure $\ov{\nabla}:E\ra E$ on $E$
the formula
$$\ov{\nabla}(e\otimes e_0)=\frac{1}{\rk(E_0)}\ov{\nabla}(e)\otimes e_0
+e\otimes\ov{\nabla}_z(e_0)$$
defines a holomorphic structure on $E\otimes_{A_{g\th}}E_0$.

\end{prop}

\Pf . (a) This follows from the explicit formulas for the action of $A_{g\th}$ on $E$.

\noindent
(b) The identities of (a) for $b=U_1$ and $b=U_2$ can be rewritten as follows:
$$\ov{\nabla}_{z-\frac{\tau}{\rk(E)}}(U_1e)=U_1\ov{\nabla}_z(e),$$
$$\ov{\nabla}_{z-\frac{1}{\rk(E)}}(U_2e)=U_2\ov{\nabla}_z(e).$$
Therefore, $U_1$ (resp. $U_2$) induces an equivalence between $\ov{\nabla}_z$ and
$\ov{\nabla}_{z-\frac{\tau}{\rk(E)}}$ (resp. $\ov{\nabla}_{z-\frac{1}{\rk(E)}}$).

Conversely, assume that $\ov{\nabla}_{z'}$ is equivalent to $\ov{\nabla}_z$.
Then $z'=z-\de_{\tau}(u)u^{-1}/\rk(E)$ for some $u\in A_{g\th}^*$.
But this implies that $\de_{\tau}(u)$ is proportional to $u$, which is possible only
for $u$ of the form $c\cdot U_1^{n_1}U_2^{n_2}$. Hence,
constant elements $z$ and $z'$ belong to the same orbit if and only if
$z-z'\in \rk(E)^{-1}(\Z+\tau\Z)$.

\noindent
(c) This follows easily from (a).
\ed

It is easy to see that not every holomorphic structure on a basic module is equivalent to
a standard one. Indeed, if $\theta$ is irrational then every basic module $E$ is isomorphic
to a direct sum $E_1\oplus E_2\oplus\ldots E_n$ of two or more other basic 
modules (this follows from Rieffel's cancellation theorem, see \cite{R-ct}).
Then the direct sum of holomorphic structures on $E_i$'s will not be equivalent
to a standard holomorphic structure on $E$, since it will admit holomorphic endomorphisms
that are not proportional to the identity, while for standard holomorphic structures
this is not the case(as follows from Corollary \ref{holhomcor} below).
Holomorphic bundles isomorphic to standard ones can be characterized as those
bundles admitting a compatible unitary connection with constant curvature (and such
that the underlying projective module is basic). In the commutative case the theorem of
Narasimhan and Seshadri \cite{NS} asserts that these are exactly stable bundles. In case of elliptic curves an equivalent characterization is that
these are bundles with only scalar endomorphisms.
One may hope that this remains true for noncommutative tori. One natural approach
to this problem would be to mimick the method of Donaldson in \cite{D}. Roughly speaking
the idea is to apply the gradient flow of the Yang-Mills functional and then use
the theorem of Connes-Rieffel \cite{CR} stating that the minima of this functional 
are exactly connections of constant curvature (this approach is advocated in \cite{Sp}). 
However, the induction
in rank employed in \cite{D} cannot be used for irrational $\theta$, since ranks
of bundles on $T_{\th}$ can be arbitrarily small. It is possible that one can modify
the argument using Rieffel's description of all critical sets of the Yang-Mills functional
(see \cite{R}).
Another possible approach would be to use
an analogue of the Fourier-Mukai transform defined in section \ref{Fouriersec}. 

\subsection{Tensor products and $\Hom$ in holomorphic category}

Now we are going to study the relation between the pairings $t_{g_1,g_2}$ and
holomorphic structures.

%so $\rk(g,\th)\cdot\ov{\nabla}$ can be considered
%as a holomorphic structure on $E_g(\th)$ considered as a left $A_{g\th}$-module (where
%$T_{g\th}$ is equipped with the complex structure $\tau$).

\begin{prop}\label{holhomprop} (a)
For every $g_1,g_2\in\SL_2(\Z)$, $z_1,z_2\in\C$ one has 
\begin{equation}\label{nablaid}
\ov{\nabla}_{z_1+z_2}(t_{g_1,g_2}(f_1\otimes f_2))=
t_{g_1,g_2}(\frac{1}{\rk(g_2,\th)}\ov{\nabla}_{\rk(g_2,\th)z_1}(f_1)\otimes f_2)+
t_{g_1,g_2}(f_1\otimes\ov{\nabla}_{z_2}(f_2)),
\end{equation}
where $f_1\in E_{g_1}(g_2\th)$, $f_2\in E_{g_2}(\th)$.

\noindent (b) The pairing $t_{g_1,g_2}$ induces a well-defined pairing
$$H^i(E_{g_1}(g_2\th),\ov{\nabla}_{\rk(g_2,\th)z_1})\otimes H^j(E_{g_2}(\th),\ov{\nabla}_{z_2})
\ra H^{i+j}(E_{g_1g_2}(\th),\ov{\nabla}_{z_1+z_2})$$
where $(i,j)$ is either $(0,0)$, $(1,0)$ or $(0,1)$.
\end{prop}

\Pf . (a) Since $\ov{\nabla}_z(f)=\ov{\nabla}_0(f)+zf$,
it suffices to consider the case $z_1=z_2=0$. If $g_1\neq 1$, $g_2\neq 1$ and
$g_1g_2\neq 1$ then (\ref{nablaid}) follows immediately from
the explicit formula for the pairing $t_{g_1,g_2}$. For $g_2=1$ it 
reduces to Leibnitz identity, while for $g_1=1$ it follows from Lemma \ref{leftLeiblem}.
Finally, in the case $g_1g_2=1$ the equality (\ref{nablaid}) follows from the identity
$$b(\frac{1}{\rk(g,\th)}\ov{\nabla}_0(f_1)\otimes f_2)+b(f_1\otimes\ov{\nabla}_0(f_2))=0$$
for the pairing (\ref{bpair}), which is easy to check using integration by parts.

\noindent
(b) This follows immediately from (a).
\ed

\begin{cor}\label{holhomcor}
Let $E=E_g^z(\th)$ and $E'=E_{g'}^{z'}(\th)$ be a pair of basic modules equipped with holomorphic
structures $\ov{\nabla}_z$ and $\ov{\nabla}_{z'}$ respectively. 

\noindent (a) The formula 
$$\ov{\nabla}(\phi)(e)=\rk(E)[\ov{\nabla}_{z'}(\phi(e))-\phi(\ov{\nabla}_z e)]$$
defines a holomorphic structure on the right $A_{g\th}$-module
$\Hom_{A_{\th}}(E,E')$ (where $\phi\in\Hom_{A_{\th}}(E,E')$, $e\in E$).

\noindent (b) Under the isomorphism (\ref{homisom}) the above holomorphic structure
on $\Hom_{A_{\th}}(E,E')$ corresponds to the operator $\ov{\nabla}_{\rk(E)(z'-z)}$
on $E_{g'g^{-1}}(g\th)$, i.e. we have
$$\Hom_{A_{\th}}(E_g^z(\th),E_{g'}^{z'}(\th))\simeq E_{g'g^{-1}}^{\rk(g,\th)(z'-z)}(g\th)$$
in a way compatible with holomorphic structures.

\noindent (c) The subspace
$H^0(\Hom_{A_{\th}}(E,E'))\sub \Hom_{A_{\th}}(E,E')$
coincides with the subspace of holomorphic morphisms $E\ra E'$.
\end{cor}

In the same way as in commutative geometry we can interpret the space
$H^1(\Hom_{A_{\th}}(E,E'))$ in terms of holomorphic extensions.

\begin{prop}\label{extprop}
Let $E=E_g^z(\th)$ and $E'=E_{g'}^{z'}(\th)$ be a pair of basic modules equipped with
standard holomorphic structures. Then there is a functorial bijection
between the space $H^1(\Hom_{A_{\th}}(E',E))$ and isomorphism classes
of extensions $0\ra E\ra F\ra E'\ra 0$ in the category of holomorphic
bundles on $T_{\th}$.
\end{prop}

\Pf . Let $0\ra E\ra F\ra E'\ra 0$ be an extension in
the holomorphic category. Since $E'$ is projective we can choose
a (non-holomorphic) splitting $s:E'\ra F$ which induces an isomorphism
of $A_{\th}$-modules $F\simeq E\oplus E'$. The holomorphic structure
$\ov{\nabla}_F$ on $F$ should have form
$\left(\matrix \ov{\nabla}_E & \phi \\ 0 & \ov{\nabla}_{E'}\endmatrix\right)$
for some $\phi\in\Hom_{A_{\th}}(E',E)$. All possible 
splittings $s:E'\ra F$ form a principal homogeneous space over
$\Hom_{A_{\th}}(E',E)$ and changing of $s$ by an element
$f\in\Hom_{A_{\th}}(E',E)$ leads to the change $\phi\mapsto\phi+
\ov{\nabla}(f)$, where $\ov{\nabla}$ is the holomorphic structure
on $\Hom_{A_{\th}}(E',E)$ defined above. This easily implies the
assertion.
\ed

\subsection{Holomorphic vectors and cohomology}

The cohomology spaces $H^*(E_g^z(\th)):=H^*(E_g(\th),\ov{\nabla}_z)$ are easy to compute.
Namely, the following result holds.

\begin{prop}\label{cohprop} 
Assume that $\Im(\tau)<0$. Let $E=E_g^z(\th)$ be a basic module equipped with
a holomorphic structure $\ov{\nabla}_z$.

\noindent
(a) If $\mu(E)>0$ then $H^1(E)=0$ and
$H^0(E)$ has dimension $\deg(E)$.

\noindent
(b) If $\mu(E)<0$ then $H^0(E)=0$ and
$H^1(E,\ov{\nabla}_z)$ has dimension $\deg(E)$.

\noindent
(c) If $z\neq 0$ then $H^*(A_{\th},\ov{\nabla}_z)=0$.

\noindent
(d) The spaces $H^0(A_{\th},\ov{\nabla}_0)$ and $H^1(A_{\th},\ov{\nabla}_0)$
are $1$-dimensional.
\end{prop}

\Pf . In the case $\deg(E)\neq 0$ the problem reduces to the computation of
the cohomology of the operator 
$$f\mapsto f'+2\pi i(\tau\mu(E)x+z)f$$
on the Schwartz space $\SS(\R)$. Now one has to use the fact that for
every complex number $a$ with $\Re(a)\neq 0$ the operator $f\mapsto f'+(ax+z)f$ on
$\SS(\R)$ has $1$-dimensional kernel and no cokernel (resp., $1$-dimensional
cokernel and no kernel) if $\Re(a)>0$ (resp. $\Re(a)<0$). To prove this fact one
can use conjugation by the operator $f\mapsto f\exp(i\Im(a)x^2/2)$ and rescaling the variable
$x$ to reduce oneself to the case $a=\pm 1$. Furthermore, making the change
of the variable $x\mapsto x\pm\Re(z)$ and conjugating by the operator
$f\mapsto f\exp(i\Im(z)x)$ we can assume that $z=0$. Now
writing $f$ as a linear combination of $(H_n(x)\exp(-x^2/2), n\ge 0)$, where
$H_n$ are Hermite polynomials, one easily proves that $\exp(-x^2/2)$ generates
the kernel for $a=1$ (resp., the cokernel for $a=-1$).

In the case $E=A_{\th}$ the assertion follows immediately from the
formula defining $\ov{\nabla}_z$.
\ed

\begin{cor} Assume that $\Im(\tau)<0$. Then
for every basic module $E$ equipped with some holomorphic structure
one has
$$\chi(E)=\dim H^i(E)-\dim H^{i+1}(E)=\deg(E),$$
where $i=0$ if $\rk(E)>0$ and $i=1$ if $\rk(E)<0$.
\end{cor}

The computation of \ref{cohprop} gives explicit bases in $H^*(E_g^z(\th))$.
In particular, the space of holomorphic vectors $H^0(E_g^z(\th))$
in a basic module $E_g^z(\th)$, such that $\deg(g)>0$, has a natural basis
$$\phi^z_{\a}(x,\b):=\e(-\tau\mu(E)x^2/2-zx)\delta_{\a}(\b),\ \a\in\Z/\deg(g)\Z,$$
where $\e(t):=\exp(2\pi i t)$, $\delta_{\a}$ is the delta-function at $\a\in\Z/\deg(g)\Z$.
The basis in $H^0(A_{\th},\ov{\nabla}_0)$ is the element $\phi^0=1\in A_{\th}$.

The following proposition is essentially equivalent to the main result of \cite{DS}.

\begin{prop}\label{DSprop} 
For every $g_1,g_2\in\SL_2(\Z)$ such that
$\deg(g_1)>0$, $\deg(g_2)>0$, $\rk(g_1,\th)>0$, $\rk(g_2,\th)>0$ one has
$\deg(g_1g_2)>0$ and
$$t_{g_1,g_2}(\phi^{\rk(g_2,\th)z_1}_{\a_1}\otimes\phi^{z_2}_{\a_2})=
\sum_{\a\in\Z/c_{12}\Z}c_{\a_1,\a_2}^{\a}\phi_{\a}^{z_1+z_2},\ \text{where}$$
$$c_{\a_1,\a_2}^{\a}=\sum_{m\in I_{g_1,g_2}(\a_1,\a_2,\a)}
\e[\frac{-\tau m^2/2+(c_1z_2-\rk(g_1g_2,\th)c_2z_1)m}{c_1c_2c_{12}}],$$
where we use the notation of Proposition \ref{pairprop}.
\end{prop}

\Pf . The fact that $\deg(g_1g_2)>0$ follows from the identity
$$\deg(g_1g_2)=\deg(g_1)\rk(g_2^{-1},\th)+\deg(g_2)\rk(g_1,\th)$$
which is the particular case of (\ref{degid}). Hence, $g_1g_2\neq 1$, 
and we can use the formula of Proposition \ref{pairprop}(a) to compute 
$t_{g_1,g_2}(\phi^{\rk(g_2,\th)z_1}_{\a_1}\otimes\phi^{z_2}_{\a_2})$. This
gives the above formula for $c_{\a_1,\a_2}^{\a}$.
\ed

The main thing to observe in the formula of the above proposition is that $\th$ enters only 
through the combination $\deg(g_1)z_2-\rk(g_1g_2,\th)\deg(g_2)z_1$.

\subsection{Translations, tensorings, pull-back and push-forward with respect to isogenies}

To conclude this section let us consider noncommutative
analogues of some natural
operations on holomorphic vector bundles.
For every $(v_1,v_2)\in\R^2$ we denote by $\a_{v_1,v_2}$ the automorphism
$U_1\mapsto\e(v_1)U_1$, $U_2\mapsto \e(v_2)U_2$ of $A_{\th}$. Then
for a right $A_{\th}$-module $E$ we denote by
$t_{(v_1,v_2)}^*E$ the space $E$ with the new right action of $A_{\th}$
given by $f\mapsto f\cdot \a_{v_1,v_2}^{-1}(a)$, where $a\in A_{\th}$.
Since $\de_{\tau}$ commutes with $A_{v_1,v_2}$, a holomorphic
structure $\ov{\nabla}:E\ra E$ induces a holomorphic structure on
$t_{(v_1,v_2)}^*E$ (given by the same operator).

Another natural operation is an analogue of tensoring with
a topologically trivial line bundle. Clearly, analogues of such line bundles 
are modules $E_1^{\la}(\th)$, where $\la\in\C$ 
(recall that $E_1^{\la}(\th)$ is $A_{\th}$ equipped with the holomorphic
structure $\ov{\nabla}_{\la}=\de_{\tau}+2\pi i\la$). It is easy to see that
for every right
$A_{\th}$-module $E$ with holomorphic structure $\ov{\nabla}$ the
tensor product $E\otimes_{A_{\th}}E_1^{\la}(\th)$ is isomorphic to $E$ as
a $A_{\th}$-module, but the holomorphic structure on it 
(defined in Proposition \ref{leftLeiblem}(c))
corresponds to $\ov{\nabla}+2\pi i\la$ under this isomorphism.
Thus, we have
$$E_g^z(\th)\otimes_{A_{\th}}E_1^{\la}\simeq E_g^{z+\la}.$$
Therefore, using Proposition \ref{leftLeiblem}(b),
we obtain an
isomorphism 
$$E\otimes_{A_{\th}} E_1^{\frac{m+n\tau}{\rk(E)}}\simeq E$$ 
for $m,n\in\Z$.

For a standard holomorphic bundle $E=E_g^z(\th)$
such that $\deg(E)\neq 0$, the map
$$f(x)\mapsto \e(\mu(E)v_1x)f(x+v_2)$$
induces an isomorphism
\begin{equation}\label{transtenseq}
t_{(v_1,v_2)}^*E \wt{\ra} E\otimes E_1^{\mu(E)(\tau v_2-v_1)}
\end{equation}
compatible with holomorphic structures.
On the other hand, we have isomorphisms
$$t_{(v_1,v_2)}^*E_1^z\simeq E_1^z$$
for all $(v_1,v_2)\in\R^2$ and all $z\in\C$ given by the map $\a_{v_1,v_2}:E_1\ra E_1$.

Note that if $c=\deg(E)\neq 0$ then taking in (\ref{transtenseq})
$(v_1,v_2)=(m/c,n/c)$ with $m,n\in\Z$ we obtain
isomorphisms
$$t_{(m/c,n/c)}^*E\wt{\ra} E$$
of holomorphic vector bundles on $T_{\th,\tau}$. Therefore, combining
induced isomorphisms of cohomology spaces with the
natural identifications $H^*(t_{(m/c,n/c)}^*E)=H^*(E)$
we obtain a collection of invertible operators
$$U_{(m/c,n/c)}:H^i(E)\ra H^i(E),$$
where $i$ is the unique degree such that $H^i(E)\neq 0$.
Similar to the commutative case we have the following result.

\begin{lem}\label{Heislem} 
The operators $U_{(m/c,n/c)}$ define an action
on $H^i(E)$ of the Heisenberg group $H$ which is a central
extension of $(\Z/c\Z)^2$ by $U(1)$, such that $U(1)$ acts
in the standard way. This representation is irreducible.
\end{lem}

\Pf . Unravelling the definitions we get the following formula for $U_{(m/c,n/c)}$
(up to a slight rescaling):
$$U_{(m/c,n/c)}f(x,\a)=\e(-\frac{m}{c}\a)f(x,\a-na),$$
where $a\in\Z$ is relatively prime to $c$. Hence, the action of these
operators on $H^i(E)$ is equivalent to the standard representation of the
Heisenberg group on the space of functions on $\Z/c\Z$.
\ed

There are also analogues of pull-back and push-forward with respect to isogenies
in our situation.
Namely, for every positive integer $N$ we have a noncommutative
morphism $\pi:T_{\th}\ra T_{N\th}$ corresponding to the embedding of algebras
$$A_{N\th}\ra A_{\th}:U_1\mapsto U_1^N, U_2\mapsto U_2.$$
If $\tau\in\C\setminus\R$ then $\pi$ can be considered as a morphism
of complex noncommutative tori $T_{\th,\tau}\ra T_{N\th,N\tau}$. Indeed,
the derivation $\de_{\tau}$ preserves the subalgebra $A_{N\th}\sub A_{\th}$
and $\de_{\tau}|_{A_{N\th}}=\de_{N\tau}$.
For a projective right $A_{\th}$-module $E$ we set $\pi_*(E)=E$, considered as
a $A_{N\th}$-module. Clearly, a holomorphic structure on $E$ induces a holomorphic
structure on $\pi_*(E)$ (with respect to $N\tau$). The analogue of pull-back
is the operation $E\mapsto \pi^*E=E\otimes_{A_{N\th}} A_{\th}$ on projective $A_{N\th}$-modules.
Here is another description of $\pi^*E$: it is the space of functions
$f:\Z\ra E$ satisfying the equation $f(n+N)=f(n)U_1^{-1}$ with the action of $A_{\th}$ given
by $fU_1(n)=f(n-1)$, $fU_2(n)=\e(n\th)f(n)U_2$. A holomorphic structure 
$\ov{\nabla}:E\ra E$ satisfying the Leibnitz rule with respect to $\de_{N\tau}$ induces
a holomorphic structure on $\pi^*E$ by the formula
$$\ov{\nabla}(f)(n)=\ov{\nabla}(f(n))+2\pi i n\tau f(n),$$
where $n\in\Z$. It is easy to check that $(\pi^*,\pi_*)$ is an adjoint pair of
functors between the categories of holomorphic vector bundles on $T_{N\th,N\tau}$ and $T_{\th,\tau}$.

\begin{prop}\label{isogprop} 
(a) For every $m$ relatively prime with $N$ one has an isomorphism
$$\pi_*E_{1,m}^z(\th)\simeq E_{N,m}^{z}(N\th)$$
of holomorphic vector bundles on $T_{N\th, N\tau}$.

\noindent (b) For every $m$ one has an isomorphism 
$$\pi^*E_{1,m}^z(N\th)\simeq E_{1,mN}^{z}(\th)$$
of holomorphic vector bundles on $T_{\th,\tau}$.
\end{prop}

\Pf . (a) The required isomorphism is given by $f(x,\a)\mapsto f(x,N\a)$.

\noindent (b) We can consider elements of $\pi^*E_{1,m}^z(N\th)$
as functions $f(x,\a,n)$ of $x\in\R$, $\a\in\Z/N\Z$ and $n\in\Z$ satisfying
$f(x,\a,n+N)=f(x+(1+mN\th)/m,\a+1,n)$. The isomorphism with $E_{1,mN}^z(\th)$
sends such a function to $g(x,n)=f(x-\frac{n(1+mN\th)}{mN},0,n)$,
where $x\in\R$ and $n\in\Z/mN\Z$.
\ed

\section{Derived categories of holomorphic bundles}

\subsection{Definition and basic properties}

Let us fix a complex number $\tau$ such that $\Im(\tau)<0$.
We are going to define certain dg-category 
$\CC(\th,\tau)$ for the complex noncommutative torus $T_{\th,\tau}$.
The corresponding cohomology category $H^0(\CC(\th,\tau))$ 
can be considered as an analogue of the derived
category of holomorphic vector bundles.

Objects of the dg-category $\CC=\CC(\th,\tau)$ are of the form
$E[n]$, where $E$ is a finitely generated projective right $A_{\th}$-module
such that $\rk(E)>0$,
equipped with a complex structure (compatible with $\de_{\tau}$),
$n$ is an integer.
We will often write $E$ instead of $E[0]$. Morphisms are defined as follows:
$$\Hom_{\CC}^{\bullet}(E[n],E'[n'])=
(\Hom_{A_{\th}}(E,E')\ra\Hom_{A_{\th}}(E,E'))[n'-n]$$
where this two-term complex is placed in degrees $n-n'$ and $n-n'+1$,
the differential sends $f\in\Hom_{A_{\th}}(E,E')$ to the $A_{\th}$-linear map
$e\mapsto\ov{\nabla}(\phi(e))-\phi(\ov{\nabla}(x))$.
The composition 
$$\Hom_{\CC}^{\bullet}(E_2,E_3)\otimes\Hom_{\CC}^{\bullet}(E_1,E_2)\ra
\Hom_{\CC}^{\bullet}(E_1,E_3)$$
is defined in the obvious way.

We denote by $\CC^{st}=\CC^{st}(\th,\tau)$ the full subcategory of $\CC$ consisting
of objects $E[n]$ such that $E$ is a standard holomorphic bundle.
Note that in this case the differential on 
$\Hom_{A_{\th}}(E,E')$ 
is proportional to the holomorphic structure $\ov{\nabla}$ on it
defined in Corollary \ref{holhomcor}.

\begin{rem} One can also consider the $\Z/2\Z$-graded analogue of $\CC$. 
Recall that the $A_{g\th}$-module $\Hom_{A_{\th}}(E,E')$
can be either even or odd (depending on the sign of its rank). We can equip
the space $\Hom_{\CC}(E,E')$ with the induced $\Z/2\Z$-grading, so that
the differential is odd.
\end{rem}

We are going to check that
in the case $\th=0$ the category $H^0\CC^{st}$ is equivalent to the full
subcategory in the derived category $D^b(X)$
of coherent sheaves on the elliptic
curve $X=X_{\tau}=\C/(\Z+\tau\Z)$. Namely, 
let us denote by $\CC^{st}(X)$ the dg-category with objects $E[n]$, where
$E$ is a stable holomorphic vector bundle on $X$, $n\in\Z$.
The spaces $\Hom_{\CC^{st}(X)}(E_1,E_2)$ are given by the Dolbeault complexes of
$E_1^{\vee}\otimes E_2$. Then this definition is extended to objects of
the form $E[n]$ as before. The associated usual category $H^0\CC^{st}(X)$
is the full subcategory in the derived category $D^b(X)$.

\begin{prop}\label{commcaseprop} 
The dg-categories $\CC^{st}(0,\tau)$ and $\CC^{st}(X_{\tau})$
are equivalent. If $\Im(\tau)<0$ then
under this equivalence $E_{n,m}^z$ corresponds to
a stable vector bundle of rank $n$ and degree $m$ on $X_{\tau}$. 
\end{prop}

\Pf . Let us define real coordinates $(x,y)$ on $\C$ by setting $z=x-\tau y$.
Then we can identify the algebra of smooth functions on the torus
$X=X_{\tau}$ with $A_0$ by setting $U_1=\exp(2\pi ix)$, $U_2=\exp(2\pi iy)$.
Then the derivation $\de_{\tau}$ on $A_0$ is proportional to $\ov{\partial}$.
Given a holomorphic vector bundle $V$ over $X$ the space $C^{\infty}(V)$ of
smooth sections of $V$ is a projective $A_0$-module. Furthermore, trivializing
the bundle of $(0,1)$-forms on $X$ we can consider the operator
$\ov{\partial}:C^{\infty}(V)\ra C^{\infty}(V)$ (appropriately rescaled)
as a holomorphic structure on the $A_0$-module $C^{\infty}(V)$. 

Let us check that if $V$ is stable then $C^{\infty}(V)$ is a basic $A_0$-module
with a standard holomorphic structure.
Assume first that $V$ is a line bundle of degree $c\neq 0$. Every such line bundle on $X$
is isomorphic to a line bundle $L_c(u)$, where $u\in\C$,
such that $C^{\infty}(L_c(u))$ is the
space of smooth functions $f$ on $\C$ satisfying the equations $f(z+1)=f(z)$,
$f(z-\tau)=\e(u-cz)f(z)$, and holomorphic structure is given by the $\ov{\partial}$-operator
with respect to the complex variable $z$. Using the real coordinates $(x,y)$ such
that $z=x-\tau y$ we can write a section of $L_c(u)$ in the form
$$f(z)=\sum_{n\in\Z}f_n(y)\e(nx),$$
where $f_n(y+1)=f_{n+c}(y)\e(u+c\tau y)$. Setting 
$$g_n(y)=f_n(y)\e(-c\tau y^2/2+(c\tau/2-u)y),$$ 
we note that $g_n(y+1)=g_{n+c}(y)$, so the
function $\phi(y,\a)=g_{\a}(y-\a/c)$ depends only on $\a\mod c\Z$.
Smoothness of $f$ implies that the functions $\phi(y,\a)$ belong
to the Schwartz spaces. It is easy to check that the map $f\mapsto \phi(y,\a)$
induces an isomorphism of $A_0$-modules $C^{\infty}(L_c(u))\simeq E_{1,c}(0)$.
Furthermore, the holomorphic structure on $L_c(u)$ corresponds to the
holomorphic structure $\ov{\nabla}_{u-c\tau/2}$ (with respect to the complex
structure $\de_{\tau}$ on $A_0$). 

Next, we observe that the construction $V\mapsto C^{\infty}(V)$ is compatible
with the operation $\pi_*$, where $\pi:\C/(d\Z+\tau\Z)\ra\C/(\Z+\tau\Z)=X$
is a natural isogeny of degree $d$.
Since every stable holomorphic vector bundle of rank $d$
has form $\pi_*L_c(u)$, applying Proposition \ref{isogprop} we derive that
the space of smooth sections of such bundle can be identified with a basic module over
$A_0$ (and that the holomorphic structure on it is standard).
\ed

Let us fix a basic right $A_{\th}$-module 
$E_0=E_g(\th)$ and equip it with a standard
holomorphic structure $\ov{\nabla}_0$.
Then as we have shown in Proposition \ref{leftLeiblem},
the map $E\mapsto E\otimes_{A_{g\th}}E_0$ extends to a functor
from the category of holomorphic bundles over $T_{g\th}$ to
the category of holomorphic bundles over $T_{\th}$.

\begin{prop}\label{holmoritaprop} 
The above functor extends to an equivalence of
dg-categories $\CC(g\th,\tau)\wt{\ra}\CC(\th,\tau)$ (resp.,
$\CC^{st}(g\th,\tau)\wt{\ra}\CC^{st}(\th,\tau)$).
\end{prop}

\Pf . Indeed, for every holomorphic vector bundles
$E$, $E'$ over $T_{g\th}$ we have a natural map
$$\Hom_{A_{g\th}}(E,E')\ra
\Hom_{A_{\th}}(E\otimes_{A_{g\th}}E_0,E'\otimes_{A_{g\th}}E_0)$$
commuting with the differentials. The inverse $A_{\th}-A_{g\th}$-bimodule
to $E_0$ will define an inverse map.
\ed

The following result is an analogue of Serre duality for the category
$H^*(\CC^{st}(\th,\tau))$.

\begin{prop}\label{Serreprop} 
(a) The pairing (\ref{bpair}) induces a perfect pairing
$$H^i(E_{g^{-1}}^{-\rk(g,\th)z}(g\th))\otimes H^{1-i}(E_g^z(\th))\ra\C$$
for $i=0,1$.

\noindent
(b) There is a canonical functorial isomorphism 
$H^i\Hom_{\CC}(E,E')\simeq H^{1-i}\Hom_{\CC}(E',E)^*$,
where $E,E'\in\CC^{st}$.
\end{prop}

\Pf . (a) The identity
$$b(\frac{1}{\rk(g,\th)}\ov{\nabla}_{-\rk(g,\th)z}(f_1)\otimes f_2)+
b(f_1\otimes\ov{\nabla}_z(f_2))=0,$$
where $f_1\in E_{g^{-1}}(g\th)$, $f_2\in E_g(\th)$, ensures that
$b$ descends to the required pairing. We already know
that dimensions of both spaces are the same. To prove that this
pairing is non-degenerate it suffices to check non-degeneracy
of the pairing
$$\ker(d/dx+x:\SS(\R)\ra\SS(\R))\otimes\coker(d/dx-x:\SS(\R)\ra\SS(\R))\ra\C$$
induced by $(f_1,f_2)=\int_{x\in\R}f_1(x)f_2(x)dx$.
But both spaces are generated by $\exp(-x^2/2)$ and $(\exp(-x^2/2),\exp(-x^2/2))\neq 0$.

\noindent (b) This follows from (a) and from the isomorphisms of Corollary
\ref{holhomcor}(b). 
\ed

\subsection{Comparison of the categories for different $\th$}

Proposition \ref{DSprop} gives explicit formulas
for the structure constants of the composition law
in the categories $H^*\CC^{st}(\th,\tau)$. Looking at these
formulas one can make the following observation.
Let us denote by $\CC^{st}(\th,\tau)_{\th'}$ the full dg-subcategory of
$\CC^{st}(\th,\tau)_{\th'}$ formed by objects $E_g^z(\th)[n]$ with
$\rk(g,\th')\neq 0$. 

\begin{thm}\label{twothetathm} 
For every pair $(\th,\th')$ the graded categories
$H^*\CC^{st}(\th,\tau)_{\th'}$ and $H^*\CC^{st}(\th',\tau)_{\th}$ are equivalent.
\end{thm}

\Pf . The main ingredient of the proof is the following computation of the composition
in the cohomology category of
$\CC^{st}(\th)=\CC^{st}(\th,\tau)$. Assume that $\deg(g_2g_1^{-1})>0$,
$\deg(g_3g_2^{-1})>0$, and $\rk(g_i,\th)>0$ for $i=1,2,3$. Then the composition
$$H^0\Hom(E_{g_2}^{z_2}(\th), E_{g_3}^{z_3}(\th))\otimes H^0\Hom(E_{g_1}^{z_1}(\th),E_{g_2}^{z_2}(\th))
\ra H^0\Hom(E_{g_1}^{z_1}(\th),E_{g_3}^{z_3}(\th))$$
can be identified using (\ref{homisom}) and Corollary \ref{holhomcor} with the map
$$H^0(E_{g_3g_2^{-1}}^{\rk(g_2,\th)(z_3-z_2)}(g_2\th))\otimes
H^0(E_{g_2g_1^{-1}}^{\rk(g_1,\th)(z_2-z_1)}(g_1\th))
\ra H^0(E_{g_3g_1^{-1}}^{\rk(g_1,\th)(z_3-z_1)}(g_1\th))$$
induced by the pairing $t_{g_3g_2^{-1},g_2g_1^{-1}}$.
Applying Proposition \ref{DSprop} we get
$$t_{g_3g_2^{-1},g_2g_1^{-1}}(\phi^{\rk(g_2,\th)(z_3-z_2)}_{\a}\otimes\phi^{\rk(g_1,\th)(z_2-z_1)}_{\b})=
\sum_{\ga}c_{\a,\b}^{\ga}\phi_{\ga}^{\rk(g_1,\th)(z_3-z_1)},\ \text{where}$$
\begin{align*}
&c_{\a,\b}^{\ga}=\sum_{m\in I_{g_3g_2^{-1},g_2g_1^{-1}}(\a,\b,\ga)}\\
&\e(\frac{-\tau m^2/2+
[\deg(g_3g_2^{-1})\rk(g_1,\th)(z_2-z_1)-\deg(g_2g_1^{-1})\rk(g_3,\th)(z_3-z_2)]m}
{\deg(g_3g_2^{-1})\deg(g_2g_1^{-1})\deg(g_3g_1^{-1})}).
\end{align*}
Therefore, $c_{\a,\b}^{\ga}$ depends on $\th$ only through the expression
$$\deg(g_3g_2^{-1})\rk(g_1,\th)(z_1-z_2)+\deg(g_2g_1^{-1})\rk(g_3,\th)(z_3-z_2).$$
Applying identity (\ref{degid}) we can rewrite this as
$$\deg(g_3g_2^{-1})\rk(g_1,\th)z_1-\deg(g_3g_1^{-1})\rk(g_2,\th)z_2+\deg(g_2g_1^{-1})\rk(g_3,\th)z_3,$$
so we see that this expression for data $(z_1,z_2,z_3,\th)$ is equal to
a similar expression for $(z'_1,z'_2,z'_3,\th')$, where $z'_i=z_i\rk(g_i,\th)/\rk(g_i,\th')$.

Now we can construct the equivalence functor
$F=F_{\th,\th'}:H^*(\CC^{st}(\th)_{\th'})\ra H^*(\CC^{st}(\th')_{\th})$.
Without loss of generality we can 
assume that $\th<\th'$. For every $g\in\SL_2(\Z)$ let us denote
$\la(g)=\rk(g,\th)/\rk(g,\th')$.
Then $F$ is defined on objects by setting 
$$F(E_g^z(\th)[n])=\cases 
E_g^{\la(g)z}(\th')[n], &\rk(g,\th')>0,\\
E_{-g}^{\la(g)z}(\th')[n-1], &\rk(g,\th')<0. \endcases$$
We have to define isomorphisms 
$$F_{g_1,g_2}:H^i\Hom_{\CC(\th)}(E^{z_1}_{g_1}(\th),E_{g_2}^{z_2}(\th))\ra 
H^{i+\eps_1-\eps_2}\Hom_{\CC(\th')}(E^{\la(g_1)z_1}_{(-1)^{\eps_1}g_1}(\th'),
E_{(-1)^{\eps_2}g_2}^{\la(g_2)z_2}(\th'))$$
compatible with the composition, where $\eps_i\in\{0,1\}$ are defined by
$(-1)^{\eps_i}=\sign\rk(g_i,\th')$. Assume first that $\deg(g_2g_1^{-1})>0$. 
Then we claim that only the following three cases can occur:
(i) $\rk(g_1,\th')>0$, $\rk(g_2,\th')>0$; 
(ii) $\rk(g_1,\th')<0$, $\rk(g_2,\th')>0$;
(iii) $\rk(g_1,\th')<0$, $\rk(g_2,\th')<0$.
Indeed, this is clear from the following geometric interpretation of the 
conditions $\rk(g,\th)>0$, $\deg(g_2g_1^{-1})>0$. 
Consider plane vectors $v_{\th}=(-\th,1)$ and $v_g=(d,c)$,
where $g=\left(\matrix a & b \\ c & d\endmatrix\right)$. 
Then $\rk(g,\th)>0$ means
that $v_g$ belongs to the half-plane $H_{\th}$ consisting of vectors $v$ such that
the pair $(v,v_{\th})$ is positively oriented. The assumption $\th<\th'$ is equivalent
to the condition $-v_{\th'}\in H_{\th}$. On the other hand,
the assumption $\deg(g_2g_1^{-1})>0$ means that the pair $(v_{g_1},v_{g_2})$
is positively oriented. Since both vectors $v_{g_1}$ and $v_{g_2}$
belong to the half-plane $H_{\th}$, our claim follows.
In the case (i) we define the map 
$$F_{g_1,g_2}:H^0(E^{\rk(g_1,\th)(z_2-z_1)}_{g_2g_1^{-1}}(g_1\th))\ra
H^0(E^{\rk(g_1,\th')(\la(g_2)z_2-\la(g_1)z_1)}_{g_2g_1^{-1}}(g_1\th'))$$
by setting 
$$F_{g_1,g_2}(\phi_{\a}^{\rk(g_1,\th)(z_2-z_1)})=\phi_{\a}^{\rk(g_1,\th')(\la(g_2)z_2-\la(g_1)z_1)}$$
for all $\a\in \Z/\deg(g_2g_1^{-1})\Z$.
In the case (ii) we define the map
$$F_{g_1,g_2}:H^0(E^{\rk(g_1,\th)(z_2-z_1)}_{g_2g_1^{-1}}(g_1\th))\ra
H^1(E^{-\rk(g_1,\th')(\la(g_2)z_2-\la(g_1)z_1)}_{-g_2g_1^{-1}}(-g_1\th'))$$
by the formula
$$F_{g_1,g_2}(\phi_{\a}^{\rk(g_1,\th)(z_2-z_1)})=
\psi_{-d_{12}\a}^{-\rk(g_1,\th')(\la(g_2)z_2-\la(g_1)z_1)}$$
where $g_2g_1^{-1}=\left(\matrix a_{12} & b_{12} \\ c_{12} & 
d_{12}\endmatrix\right)$,
$(\psi_{\a}^{-\rk(g,\th)z},\a\in\Z/c\Z)$ denotes the basis of 
$H^1(E_{g^{-1}}^{-\rk(g,\th)z}(g\th))$
dual to the basis $(\phi_{\a}^z)$ of $H^0(E_g^z(\th))$
with respect to the natural pairing (see \ref{Serreprop}).
In the case (iii) the definition is similar to the case (i).
We can easily extend these definitions to arbitrary morphisms using the compatibility
with Serre duality.

It remains to check the compatibility of this functor with the composition
$$H^j\Hom_{\CC(\th)}(E^{z_2}_{g_2}(\th),E_{g_3}^{z_3}(\th))\otimes
H^i\Hom_{\CC(\th)}(E^{z_1}_{g_1}(\th),E_{g_2}^{z_2}(\th))\ra
H^{i+j}\Hom_{\CC(\th)}(E^{z_1}_{g_1}(\th),E_{g_3}^{z_3}(\th))$$
and the composition of the corresponding objects in $\CC^{st}(\th')$.
By Serre duality it suffices to consider the case when
$\deg(g_2g_1^{-1})>0$, $\deg(g_3g_2^{-1})>0$ and either $\rk(g_i,\th')>0$ for
$i=1,2,3$, or $\rk(g_1,\th')<0$, $\rk(g_2,\th')>0$, $\rk(g_3,\th')>0$.
In the former case this compatibility follows from our observation in the beginning
of the proof. In the latter case we have to use some cyclic symmetry of the 
structure constants $c_{\a,\b}^{\ga}$ appearing above. To reflect the dependence
of these constants on various data let us change the notation to
$c_{\a,\b}^{\ga}(g_1,g_2,g_3;z_1,z_2,z_3;\th)$. The required compatibility follows from the
equality
$$c_{\a,\b}^{\ga}(g_1,g_2,g_3;z_1,z_2,z_3;\th)=
c_{-d_{13}\ga,\a}^{-d_{12}\b}(g_2,g_3,-g_1;\la(g_2)z_2,\la(g_3)z_3,
\la(g_1)z_1;\th'),$$
where $g_jg_i^{-1}=\left(\matrix a_{ij} & b_{ij} \\ c_{ij} & d_{ij}\endmatrix\right)$.
This equality can be easily checked using Lemma \ref{pairlem}(c).
\ed

\begin{cor}
If $\th$ and $\th'$ are irrational then the graded categories
$H^*(\CC^{st}(\th,\tau))$ and $H^*(\CC^{st}(\th',\tau))$ are equivalent.
\end{cor}

In the case when $\th$ is irrational
and $\th'=0$ we can extend the composition of the above functor
$H^0\CC^{st}(\th)_{\th'}\ra H^0\CC^{st}(0)$ with the embedding $H^0\CC^{st}(0)\sub D^b(X_{\tau})$
(constructed in Proposition \ref{commcaseprop})
to the entire category $H^0\CC^{st}(\th)$.
Namely, assuming that $\th<0$ we extend  
the functor $F$ above to $E_g^z(\th)$ where
$g=\left(\matrix 0 & 1 \\ -1 & 0\endmatrix\right)$ 
by setting $F(E_g^z(\th))=\OO_{-\th z}$ (this is the structure sheaf of the
point $-\th z\mod(\Z+\tau\Z)\in X_{\tau}$).
The above construction of $F$ on morphisms then can be modified appropriately 
to prove the equivalence of
the category $H^*(\CC^{st}(\th,\tau))$ with some full subcategory of
the derived category of coherent sheaves on the elliptic curve
$\C/\Z+\tau\Z$. Below we will give a more natural construction of
the composition of this equivalence with the Fourier-Mukai transform.

\subsection{Analogue of the Fourier-Mukai transform}\label{Fouriersec}

As before we fix a complex structure $\tau$ on $T_{\th}$ such that $\Im(\tau)<0$.
Leibnitz rule implies that
for every holomorphic structure $\ov{\nabla}$ on
a right $A_{\th}$-module $E$ the cohomology spaces
$H^*(E,\ov{\nabla}+2\pi i z\id)$, where $z\in\C$, depend
only on $z\mod(\Z+\tau\Z)$. So one can try to define an analogue
of the Fourier-Mukai transform of $(E,\ov{\nabla})$ by 
descending this family of spaces to
the elliptic curve $X=X_{\tau}=\C/(\Z+\tau\Z)$.

For every $A_{\th}$-module $E$
let us consider a sheaf of $\OO$-modules
$E_{\C}$ on $\C$ such that sections of $E_{\C}$ over
an open set $U\sub\C$ are $E$-valued
holomorphic functions on $U$. We can define two morphisms
$$\rho_v:t^*_{v}E_{\C}\wt{\ra}E_{\C}$$
where $v=1$ or $\tau$; $t_v:\C\ra\C$ denotes the translation $z\mapsto z+v$.
Namely, for a local section $f$ of $E_{\C}$ we set
$$\rho_v(f)(z)=f(z+v)U^v,$$ 
where $U^v=U_1$ or $U_2$
for $v=\tau$ or $1$ respectively.
Clearly, we have
$$\rho_1\circ t^*_1\rho_{\tau}=\exp(2\pi i\th)\rho_{\tau}\circ
t^*_{\tau}\rho_1.$$
Changing $\rho_{\tau}$ to
$$\rho'_{\tau}=\exp(-2\pi i\th z)\rho_{\tau}$$
we will obtain
$$\rho_1\circ t^*_1\rho'_{\tau}=\rho'_{\tau}\circ t^*_{\tau}\rho_1.$$
Let us denote by $E_{X}$ the sheaf of
$\OO$-modules on the elliptic curve $X$ (with respect to the classical
topology) obtained from $E_{\C}$
using the $\Z^2$-action given by $\rho_1$ and $\rho'_{\tau}$. 

Now assume in addition that $E$ is equipped with a holomorphic
structure $\ov{\nabla}$. Then we define $\SS(E)=\SS(E,\ov{\nabla})$ to be
the complex 
$$E_{X}\lrar{d_{\ov{\nabla}}} E_{X}$$
of $\OO$-modules over $X$, where $d_{\ov{\nabla}}$ is defined
in terms of the corresponding operator on $E_{\C}$:
$$d_{\ov{\nabla}}(f)(z)=
\ov{\nabla}(f(z))+2\pi i zf(z).$$
We leave for the reader to check that this operator is
compatible with the $\Z^2$-action so it descends to $E_X$.

The following basic property is one of the reasons to consider $\SS$
as an analogue of the Fourier transform.

\begin{lem}\label{transtensprop}
For every $\la\in\C$, $(v_1,v_2)\in\R^2$ one has
$$\SS(E\otimes E_1^{\la})\simeq t_{\la}^*\SS(E),$$
$$\SS(t_{(v_1,v_2)}^*E)\simeq \PP_{(v_1,v_2)}\otimes\SS(E),$$
where $\PP_{(v_1,v_2)}$ is the holomorphic bundle on $X$ corresponding
to the following action of $\Z+\tau\Z$ on the trivial bundle $\OO_{\C}$:
$$\rho_{\tau}(f)(z)=\e(-v_1)f(z+\tau),\ \rho_1(f)(z)=\e(-v_2)f(z+1).$$ 
\end{lem}

The proof is the direct application of the definitions. 

Now let us compute the Fourier transform of standard objects.

\begin{prop}\label{Fourierprop}
(a) Let $E_{n,m}^0(\th)$ be a basic module corresponding
to a pair $(n,m)$ such that $m>0$ (and as usual $m\th+n>0$)
equipped with the complex structure $\ov{\nabla}_0$.
Then the complex $\SS(E_{n,m}^0(\th))$
is quasi-isomorphic to the vector bundle $V_{m,-n}$
on $X$ obtained by the descent from the vector bundle $\OO_{\C}^m$
with respect to the following action of $\Z+\tau\Z$:
$$\rho_{\tau}(f)(z,\a)=f(z+\tau,\a-1)\e(\frac{n}{m}z+\frac{\tau}{2\mu}),\
\rho_1(f)(z,\a)=f(z+1,\a)\e(-\frac{n}{m}\a),$$
where we think about sections of $\OO_{\C}^m$ as functions $f(z,\a)$
in $z\in\C$ and $\a\in\Z/m\Z$, which are holomorphic in $z$. 
The vector bundle $V_{m,-n}$ has degree $-n$ and is stable.

\noindent (b) If $m<0$ (and $m\th+n>0$)
then $\SS(E_{n,m}^0(\th))$ is quasi-isomorphic
to $V[-1]$, where $V$ is a stable vector bundle on $X$.

\noindent (c) The complex $\SS(E_{1,0}^z(\th))$ is quasi-isomorphic
to $\OO_{-z}[-1]$, where $\OO_{-z}$ is the structure sheaf  
of the point $-z\mod(\Z+\tau\Z)$ in $X$.
\end{prop}

\Pf . (a) Set $E=E_{n,m}^0(\th)$. A version of Proposition \ref{cohprop} shows that
the morphism $d_{\ov{\nabla}_0}:E_{\C}\ra E_{\C}$ is surjective
and its kernel is a vector bundle on $\C$. Moreover, we have
an explicit isomorphism of the kernel with the rank-$m$ trivial bundle on $\C$:
for a holomorphic function $f(z,\a)$ on $U\times\Z/m\Z$, where $U\subset\C$
is an open subset, the function $f(z,\a)\e(-\tau\mu(E)x^2/2-zx)$
defines a section of $\ker(d_{\ov{\nabla}_0})$ over $U$.
It is easy to check that the $\Z^2$-action on $\ker(d_{\ov{\nabla}_0})$
induced by the action of $U_1$ and $U_2$ coincides with the action
defining $V_{m,-n}$. The fact that the vector bundle $V_{m,-n}$
has degree $-n$ and is stable is well-known (note that the sign with $nz/m$
in the definition of $\rho_{\tau}$ is opposite to the standard one
due to our assumption $\Im(\tau)<0$).

\noindent (b) The proof of Proposition \ref{cohprop} shows that
the morphism $d_{\ov{\nabla}_0}:E_{\C}\ra E_{\C}$ has no kernel and
its cokernel is a bundle $V$ of rank $-m$ on $\C$. 
Let $g\in\SL_2(\Z)$ be a matrix with $(m,n)$ as the second row.
Proposition \ref{Serreprop} gives a family of isomorphisms 
$$V_z^*\simeq H^0(E_{g^{-1}}^{-\rk(g,\th)z}(g\th))$$
The identifications $V_z\simeq V_{z+1}$ and $V_z\simeq V_{z+\tau}$ correspond
to similar isomorphisms of the $H^0$-spaces induced by the left
action of generators of $A_{\th}$ on $E_{g^{-1}}(g\th)$.
In other words, the bundle dual to $V$ (or rather its descent
to $X$) is isomorphic to the Fourier-Mukai transform of $E_{g^{-1}}^0(g\th)$
considered as a left $A_{\th}$-module. Therefore, our claim follows from
the analogue of part (a) for the Fourier-Mukai transform of left $A_{\th}$-modules.

\noindent
(c) The formula for the operator
$\ov{\nabla}_z:A_{\th}\ra A_{\th}$ shows that the complex $\SS(E_{1,0}^{z_0}(\th))$
is exact outside the point $-z_0\mod(\Z+\tau\Z)\in X$ and that locally near $-z_0$
it is quasiisomorphic to the complex $\OO\stackrel{z_0+z}{\ra}\OO$.
\ed

%Remark about additional structure on complexes $\SS(E_g^0(\th))$ coming from 
%the left action of $A_{g\th}$.

\subsection{Equivalences}

Let us call an object $E\in D^b(X)$ {\it stable} if $\Hom(E,E)=\C$.
Equivalently, $E\simeq V[n]$ for some $n\in\Z$, where $V$ is either
a stable vector bundle or a structure sheaf of a point.

\begin{thm}\label{Fourierthm} Assume that $\th$ is irrational.
Then the functor $E\mapsto \SS(E)$ extends to an equivalence
of $H^0\CC^{st}(\th,\tau)$ with the full subcategory of $D^b(X)$ consisting
of stable objects.
\end{thm}

\Pf . First of all, we extend $\SS$ to all objects of $\CC^{st}(\th,\tau)$
by requiring that it commutes with the shift $E\mapsto E[1]$.
It is clear that a holomorphic map $E\ra E'$ induces a map of
complexes $\SS(E)\ra\SS(E')$. On the other hand,
a morphism $E\ra E'[1]$ in $H^0\CC^{st}(\th,\tau)$ can be interpreted
as an extension class $0\ra E'\ra F\ra E\ra 0$ (see Proposition
\ref{extprop}). Such an extension induces an exact sequence of
complexes $0\ra\SS(E')\ra\SS(F)\ra\SS(E)\ra 0$, hence a morphism
$\SS(E)\ra\SS(E')[1]$ in the derived category.
Clearly, these maps define a functor from $H^0\CC^{st}(\th,\tau)$ to
$D^b(X)$. Using Proposition \ref{Fourierprop} we see
that $\SS(E)$ is a stable object for every $E\in\CC^{st}(\th,\tau)$.
%$\dim\Hom^i(E,E')=\dim\Hom^i(\SS(E),\SS(E'))$ for
%every $E,E'$.
Let us check
that the natural maps
\begin{equation}\label{maps}
\Hom^i(E,E')\ra\Hom^i_{D^b(X)}(\SS(E),\SS(E'))
\end{equation}
are
isomorphisms. The idea of the proof is that both these $\Hom$-spaces
are irreducible representations of the same
Heisenberg group and that the map between them is a map
of Heisenberg modules. Recall that the Heisenberg group action
on $\Hom^i(E,E')=H^i(\Hom_{A_{\th}}(E,E'))$ was constructed in
Lemma \ref{Heislem}. On the other hand, the Heisenberg group
action on $\Hom^i(A,B)$ for every pair of stable objects 
in $D^b(X)$ (such that $A\not\simeq B[n]$ and $\Hom^i(A,B)\neq 0$)
is constructed as follows. The functors of translation by points of $X$
and of tensoring by topologically trivial line bundles on $X$ generate
the action of the Heisenberg groupoid $\bH$ on $D^b(X)$ which
is an extension of $X\times X$ by the groupoid of $\G_m$-torsors
(see \cite{P}). Now for every stable object $A$ in $D^b(X)$ let
$K_A\sub X\times X$ be the subgroup of translations and tensorings
preserving $A$ up to an isomorphism. Then the abelian
group $K_A\cap K_B$ is finite and has a natural Heisenberg extension
acting on $\Hom^i(A,B)$. Indeed, applying an automorphism of $D^b(X)$
we can assume that $A=\OO_X$, in which case this construction is
standard (note that $K_{\OO_X}=X\subset X\times X$ is the subgroup of all
translations). The compatibility of $\SS$ with Heisenberg group
actions on $\Hom$'s now follows easily from Proposition \ref{transtensprop}.
It remains to check that all the maps (\ref{maps}) are non-zero.
The commutative diagram
\begin{equation}
\begin{array}{ccc}
\Hom^1(E',E)\otimes\Hom(E,E') &\lrar{}&\Hom^1(E,E)\\
\ldar{}         &&\ldar{}\\
\Hom^1(\SS(E'),\SS(E))\otimes\Hom(\SS(E),\SS(E')) &\lrar{}&
\Hom^1(\SS(E),\SS(E))
\end{array}
\end{equation}
whose rows are perfect pairings, shows that it is enough to prove
nonvanishing of the maps
$$\Hom^1_{\CC^{st}(\th,\tau)}(E,E)\ra\Hom^1_{D^b(X)}(\SS(E),\SS(E)).$$
This can be done by the direct computation. Let us assume that
$E=E_{n,m}^0(\th)$, where $(n,m)\neq (1,0)$ (we leave the case $E=A_{\th}$ to the reader).
Recall that according to Proposition \ref{extprop} 
a generator of $\Hom^1(E,E)$ gives rise to a holomorphic extension
$$0\ra E\ra E^{(2)}\ra E\ra 0.$$
It is easy to see that $E^{(2)}=E\oplus E$ with the holomorphic structure
given by
$$\left(\matrix \ov{\nabla}_0 & 1\\ 0 &\ov{\nabla}_0\endmatrix\right).$$
Let us assume that $m=\deg(E)>0$ (the case $\deg(E)<0$ is considered similarly). 
Then $\SS(E)$ is the bundle $V_{m,-n}$
described in Proposition \ref{Fourierprop}(a).  
Local sections of $\SS(E^{(2)})$ over an open subset of $\C$ are
described by $E\oplus E$-valued functions of the form
$$((f_1(z,\a)+f_2(z,\a)x)\e(-\tau\mu(E)x^2/2-zx),f_2(z,\a)\e(-\tau\mu(E)x^2/2-zx)),
$$
where $f_1(z,\a)$ and $f_2(z,\a)$ are holomorphic in $z$.
The $\Z^2$-action on pairs $(f_1(z,\a),f_2(z,\a))$ is given by the rule
\begin{align*}
&\rho_{\tau}(f_1,f_2)=(f_1(z+\tau,\a-1)-\mu f_2(z+\tau,\a-1),f_2(z+\tau,\a-1))
\e(\frac{n}{m}z+\frac{\tau}{2\mu}),\\
&\rho_1(f_1,f_2)=(f_1(z+1,\a),f_2(z+1,\a))\e(-\frac{n}{m}\a),
\end{align*}
where $\mu=\mu(E)$.
This easily implies that $\SS(E^{(2)})$ is a non-split extension of $V_{m,-n}$ by itself.
\ed

\begin{rem} It is not difficult to see
that the restriction of the functor $\SS$ to $H^0\CC^{st}(\th,\tau)_0$ 
is the composition of the functor 
$F_{\th,0}$ from Theorem \ref{twothetathm}
with the embedding of $H^0\CC^{st}(0,\tau)$ in $D^b(X)$
given in Proposition \ref{commcaseprop}, followed
by the standard Fourier-Mukai transform
on $D^b(X)$ (where the elliptic curve $X$ is identified with its dual).
\end{rem}

Now let us explain the relation between the constructed equivalences
and the Morita equivalences.
Assume that $\th<0$, $\th'=g\th<0$ and $\rk(g,\th)>0$. Then
the Morita functor from $M_{g,\th}:\CC(g\th,\tau)\ra\CC(\th,\tau)$ is
given by $E\mapsto E\otimes_{A_{g\th}}E_g^0(\th)$. Hence,
it maps $E^z_h(\th')$ to $E^u_{hg}(\th)$, where $u=z/\rk(g,\th)$.
On the other hand,
the functor $F_{\th,0}$ sends an object
$E^u_{hg}(\th)$ to $E^{u\rk(hg,\th)/\rk(hg,0)}_{hg}(0)$
or $E^{u\rk(hg,\th)/\rk(hg,0)}_{-hg}(0)[1]$ depending on whether
$\rk(hg,0)$ is positive or negative. Similarly,
$F_{\th',0}(E^z_h(\th'))$ is either $E^{z\rk(h,\th')/\rk(h,0)}_h(0)$
or $E^{z\rk(h,\th')/\rk(h,0)}_{-h}(0)[1]$. Note that
$u\rk(hg,\th)=z\rk(h,\th')$. It follows that the autoequivalence 
$G=F_{\th,0}\circ M_{g,\th}\circ F_{\th',0}^{-1}$ of the subcategory
of stable objects $D^b(X)$ sends
$E^t_{h}(0)$ to $E^{t\rk(h,0)/\rk(hg,0)}_{hg}(0)$ or
to $E^{t\rk(h,0)/\rk(hg,0)}_{-hg}(0)[1]$ depending on the sign of $\rk(hg,0)$.
The autoequivalence $G$ is induced by some autoequivalence of
the derived category $D^b(X)$ that acts as the matrix $g^t$ on the column
vector $(\deg, \rk)$. 

\subsection{Tilted $t$-structures}

Let $\CC$ be an abelian category.
Recall that a torsion pair in $\CC$ is a pair of full subcategories
stable under extensions $(\CC_1,\CC_2)$, such that
for every $A_1\in\CC_1$, $A_2\in\CC_2$ one has $\Hom_{\CC}(A_1,A_2)=0$
and every object $A\in\CC$ fits into an exact sequence
$$0\ra A_1\ra A\ra A_2\ra 0$$
with $A_1\in\CC_1$, $A_2\in\CC_2$.
In this situation the above exact sequence is unique
up to a unique isomorphism (hence, it depends functorially on $A$).

A torsion pair $p=(\CC_1,\CC_2)$ gives rise to the following $t$-structure
in the derived category $D(\CC)$:
$$D^{p,\le 0}:=\{K\in D(\CC):\ H^{>0}(K)=0,\ H^0(K)\in\CC_1\},$$ 
$$D^{p,\ge 1}:=\{K\in D(\CC):\ H^{<0}(K)=0,\ H^0(K)\in\CC_2\}.$$ 
The heart of this $t$-structure $\CC^p:=D^{p,\le 0}\cap D^{p,\ge 0}$
is equipped with the torsion pair $(\CC_2[1],\CC_1)$.

In the case when $\CC$ is the category $\Coh(X)$ of coherent sheaves
on an algebraic curve $X$ over a field $k$, there are natural torsion
pairs $(\Coh_{>\th}(X),\Coh_{\le\th}(X))$ and
$(\Coh_{\ge\th}(X),\Coh_{<\th}(X))$ in $\Coh(X)$ associated with every
real number $\th$ (for irrational $\th$ these two pairs are the same). 
Namely, for a subset $I\sub\R\cup\{+\infty\}$ we denote
by $\Coh_I(X)$ the smallest full subcategory in $\Coh(X)$ closed
under extensions, containing all simple sheaves with slope $\mu\in I$
(where the slope of a stable bundle is defined as the ratio $\deg/\rk$,
a slope of a torsion sheaf is defined to be $+\infty$).

Now let us assume that $\th$ is irrational and denote by
$(D^{\th,\le 0},D^{\th,\ge 1})$ the $t$-structure on $D^b(X)$ associated
with the torsion pair $(\Coh_{>\th}(X),\Coh_{<\th}(X))$.
Let $\Coh(X)^{\th}$ be the heart of this $t$-structure, so that
$\Coh(X)^{\th}$ is equipped with the torsion pair 
$(\Coh_{<\th}[1],\Coh_{>\th})$.

\begin{prop} Assume that $\th<0$ then the functor 
$F=F_{\th,0}:H^0\CC^{st}(\th)\ra D^b(X)$ defined for the pair $(\th,0)$
sends the subcategory of holomorphic vector bundles on $T_{\th}$
to $\Coh(X)^{-\th^{-1}}[-1]$.
\end{prop}

\Pf . Indeed, if $\rk(g,0)>0$ then $F(E_g^z(\th))$ is a vector bundle
of degree $m$ and rank $n$ such that $m\th+n>0$, hence
$F(E_g^z(\th))\in\Coh_{<-\th^{-1}}$. Similarly, if $\rk(g,0)\le 0$ then 
$F(E_g^z(\th))\in\Coh_{>-\th^{-1}}[-1]$.
\ed

\end{document}